\documentclass[10pt]{amsart}
\usepackage{amsmath}
\usepackage{amssymb}
\usepackage{graphicx}
\usepackage{mathrsfs}  
\usepackage{hyperref}

\usepackage[usenames,dvipsnames]{xcolor}

\def\P{\mathbb{P}}
\def\E{\mathbb{E}}
\def\var{\operatorname{var}}
\def\beq*{\begin{eqnarray*}}
\def\eeq*{\end{eqnarray*}}

\def\be{\begin{equation}}
\def\ee{\end{equation}}
\def\eps{\varepsilon}
\def\b{\beta}

\def\a{\alpha}
\def\g{\gamma}

\newcommand\fpf{\hfill{$\blacksquare$}}

\def\NB{\operatorname{NB}}
\def\Bin{\operatorname{Bin}}
\def\N{\operatorname{N}}
\def\Pois{\operatorname{Pois}}
\newtheorem{theorem}{Theorem}

\newtheorem{remark}{Remark}

\numberwithin{subcase}{case}

\author[Pawe{\l} Hitczenko]{Pawe{\l} Hitczenko
}
\address{Department of Mathematics, Drexel University, Philadelphia, 
PA  19104, USA} 
\email{pawel.hitczenko@drexel.edu}

\title[]{On limiting distributions of Graham, Knuth, Patashnik recurrences
}

\keywords{Generating polynomial, recurrence, limiting distribution}
\subjclass[2020]{05A16, 60C05, 60F05}

\begin{document}
\maketitle
\begin{abstract}
Graham, Knuth and Patashnik in their book \emph{Concrete Mathematics} called for development of a general theory of the solutions of recurrences defined by
\[\left|{ n\atop k}\right|=(\a n+\b k+\g)\left|{n-1\atop k}\right|+(\a' n+\b' k+\g ')\left|{n-1\atop k-1}\right|+I_{n=k=0}\]
for $0\le k\le n$ and six parameters $\a,\b,\g,\a'\b',\g'$. Since then, a number of authors investigated various properties of the solutions of these recurrences.
In this note we consider a probabilistic aspect, namely we consider the limiting distributions of sequences of integer valued random variables naturally associated with the solutions of such recurrences. We will give a complete description of the limiting behavior when $\a'=0$ and the remaining five parameters are non--negative.
\end{abstract}

\section{Introduction and Motivation}\label{sec:mot}

Research problem {\bf 6.94} \cite[p.~319]{gkp} (numbered {\bf 6.89} in the first edition) calls for  developing a general theory of the solutions to  the recurrences of the  form 
\be\label{wrec} \left|{ n\atop k}\right|=(\a n+\b k+\g)\left|{n-1\atop k}\right|+(\a' n+\b' k+\g ')\left|{n-1\atop k-1}\right|+I_{n=k=0}\ee
for $0\le k\le n$ and six parameters $\a,\b,\g,\a'\b',\g'$.
Since then, a number of authors studied different aspects of the solutions of these recurrences, often under some restrictions on the parameters, see e.g.  
 \cite{bsv,m,s,ss,w}.  
Without specifically referring to \cite{gkp} recurrences of that type were also studied in \cite{n} and \cite{rr} in a fairly general setting. (Of course, for specific values
of the parameters recurrences like \eqref{wrec} were analyzed for a long time since the vast majority of classical triangles of numbers including binomial, Stirling, Eulerian, Lah, and many others, are of that form.) 

Building on an earlier work \cite{w}, the authors of \cite{bsv} classified the partial differential equations (PDEs) satisfied by the bivariate exponential generating functions (which we abbreviate to BGF) associated with \eqref{wrec} and used the method of characteristics to solve them (the solutions in most cases are expressions involving a function whose inverse is given in series expansion). Methods based  on analyzing BGFs and the  PDEs  are exploited in \eqref{m} to derive numerous results and properties of the solutions of \eqref{wrec} for various sets of values of parameters. 
 
In \cite{n}, under the assumption $\a'=0$ some explicit expressions for the solutions of \eqref{wrec} were obtained through factorization properties.
 For example,  \cite[Theorem~18]{n} can be deciphered as
\[\left|n\atop k\right|=\prod_{m=0}^{k-1}(\b'i+\g')\sum_{i=0}^n\sum_{j=0}^n\left[n\atop j\right]{j\choose i}\left\{j\atop k\right\}\a^{n-j}\g^{j-i}\b^{i-k}.\]
We note that this formula was also (but later) given by Wilf \cite[Section~2.2]{w} who noted the fact that it is factored as two polynomials, one in $\b'$ and $\g'$ only and the other in the remaining three parameters.  Under the additional assumption that $\b'=0$ and the remaining four parameters have integer values similar formulas were  also derived in \cite{rr}. Further results along these lines were obtained by Spivey~\cite{s}.
In another direction, when $\b'=0$ and the remaining five parameters are non--negative, Wilf \cite {w} showed that the generating polynomials
\be\label{poly}
P_n(x)=\sum_{k=0}^m\left|{ n\atop k}\right|x^k,\quad n\ge0.
\ee
have only real roots.

In this paper we will be interested in the asymptotic distributions of random variables naturally associated with recurrence \eqref{wrec}. Specifically, if $P_n(x)$ is given by \eqref{poly} then 
\[\frac{P_n(x)}{P_n(1)}=\sum_{k=0}^n\frac{\left|n\atop k\right|}{\sum_{j=0}^n\left|n\atop j\right|}x^k=\E\, x^{X_n},
\] 
where $X_n $ is an integer valued random variable defined by 
\be\label{Xn}\P(X_n=k)
= \frac{\left|n\atop k\right|}{\sum_{j=0}^n\left|n\atop j\right|},\quad 0\le k\le n.\ee
Our  focus will be on the limiting distribution of (suitably normalized) sequence $(X_n)$ as $n\to\infty$.
 One common means of studying the properties of the sequence $(X_n)$ is through a BGF associated to it. That is, if we let 
\be\label{F}F(z,x):=\sum_{n=0}^\infty\frac{z^n}{n!}P_n(x),\ee
then the probability generating function (PGF) of a random variable $X_n$ defined by \eqref{Xn} is 
\[p_n(x)=\E\ x^{X_n}=\frac{[z^n]F(z,x)}{[z^n]F(z,1)},\]
where as usual $[z^n](\ \cdot\ )$ denotes the coefficient of $z^n$ in  the series expansion of the expression $(\ \cdot\  )$. From now on, without mentioning it explicitly we will assume that  $(X_n)$ are related to $F(z,x)$ by the above equality. The expected value and the variance of $X_n$ will be denoted by $\mu_n$ and $\sigma_n^2$, respectively. 

Throughout the paper we will use symbols \lq\lq\ $\stackrel d\rightarrow$\ \rq\rq and \lq\lq\ $\stackrel d=$\ \rq\rq to denote, respectively,  the convergence and equality in distribution of random variables. We will denote by $\N(\nu,\sigma^2)$, $\Bin(n,p)$ and $\Pois(\lambda)$, respectively, a normal, a binomial and a Poisson distributions with the indicated parameters.

It follows from  \eqref{wrec} that the generating  polynomials \eqref{poly}
 satisfy a  recurrence 
\be\label{wpol}P_n(x)=((\a'x+\a)n+(\b'+\g')x+\g)P_{n-1}(x)+x(\b+\b'x)P'_{n-1}(x).\ee
 This is a special case of a more general scheme, namely,
 \[
 P_n(x)
=a_n(x)P_{n-1}(x)+b_n(x)P_{n-1}^{'}(x)
\]
where $(a_n)$, $(b_n)$ are sequences of known functions (typically  polynomials of low degrees)  and  $P_0(x)$ is given. Various aspects of such recurrences  have been studied throughout the years, including some recent work.
  Closely related to the subject of this paper, the authors of \cite{hcd} 
  developed the limiting theory of general Euclidean recurrences when the coefficient functions $a_n(x)$, $b_n(x)$ are  of the form:
\[a_n(x)=\a(x) n+\gamma(x),\quad b_n(x)=\beta(x)(1-x)\]
under additional conditions on $\a(x)$, $\g(x)$ and $\b(x)$.
 This corresponds to 
\[\a(x)=\a'x+\a,\quad \g(x)=(\b'+\g')x+\g,\quad \b(x)=x\b,\quad \b'=-\b\] 
and covers some of the cases.  
 In this paper we will take the next step in describing the limiting distribution of $(X_n)$. While we will not do it for a full range of all six parameters we will give a complete picture when $\a'=0$ and the remaining five parameters are non--negative.  What happens in the other cases, is an open problem.

\section{A PDE associated with \eqref{wpol}}

Differentiating \eqref{F} with respect to $z$ and using \eqref{wpol} yields
\be\label{pde}
(1-\a z-\a'xz)\frac{\partial F(z,x)}{\partial z}-x(\b+\b'x)\frac{\partial F(z,x)}{\partial x}=((\a'+\b'+\g')x+\a+\g)F(z,x).
\ee 
with $F(0,x)=1$.

In many cases this PDE  may be solved explicitly  by the method of characteristics. The first step is to  set
\be\label{char}\frac
{dz}{1-\a z-\a'xz}=-\frac{dx}{x(\b'x+\b)}=\frac{dF}{((\a'+\b'+\g')x+\a+\g)F}
\ee
which, upon solving the first of these equations 
 reduces \eqref{pde}  to the first order linear ordinary differential equation that can be explicitly solved provided the involved integrals have explicit expressions.

\section{The case $\a'=0$.}
A particularly tractable  situation arises when  $\a'=0$ in which case the first characteristic equation in \eqref{char} is separable and all integrals can be evaluated. By following the usual procedure (see e.g. Section~3.1 in \cite{hcd} for a good description) and 
 using $F(0,x)=1$ one is eventually  led to 
\be\label{Fap0}F(z,x)=\frac{(1-\a z)^{-\g/\a}}
{\left(1+\frac{\b'x}\b\left(1-(1-\a z)^{-\b/\a}\right)\right)^{1+\g'/\b'}}
\ee
 as has been derived in \cite{w} with the aid of Maple.

Having the explicit form of $F(z,x)$ makes it amenable to the asymptotic analysis. In particular, for several sets of the values of the parameters,  the algebraic singularity schema as described in \cite[Theorem~IX.12,~p.~676]{fs} applies.  Let us recall that statement.
\begin{theorem} \label{ass} {\rm (Algebraic singularity schema).}
Let $F(z,x)$ be a function that is bivariate analytic at $(z,x)=(0,0)$ and has non--negative coefficients. Assume the following conditions:
\begin{itemize}
\item[(i)] {\rm Analytic perturbation:}  there exist three functions $A$, $B$, $C$ analytic in the domain $\mathcal D=\{|z|\le r\}\times\{|x-1|<\epsilon\}$, such that, for some $r_0$ with $0<r_0\le r$ the following representation holds, with $\kappa\notin \mathbb Z_{\le0}$,
\[F(z,x)=A(z,x)+B(z,x)C(u,x)^{-\kappa};\]
furthermore, assume that, in $|z|\le r$, there exists a unique root $\rho$ of the equation $C(z,1)=0$, that this root is simple, and that $B(\rho,1)\ne0$.
\item[(ii)] {\rm Non--degeneracy:} one has $\partial_zC(\rho,1)\cdot\partial_xC(\rho,1)\ne0$, ensuring the existence of a non--constant $\rho(x)$ analytic at  $x=1$, such that $C(\rho(x),x)=0$ and $\rho(1)=\rho$.
\item[(iii)] {\rm Variability:} one has
\[\mathfrak v\left(\frac{\rho(1)}{\rho(x)}\right)\ne0,\]
where
\[\mathfrak v(f):=\frac{f''(1)}{f(1)}+\frac{f'(1)}{f(1)}-\left(\frac{f'(1)}{f(1)}\right)^2.
\]
\end{itemize}
Then, the random variable with probability generation function
\[p_n(x)=\frac{[z^n]F(z,x)}{[z^n]F(z,1)}\]
converges in distribution to a Gaussian random variable with a speed of convergence that is $O(n^{-1/2})$. The mean $\mu_n$ and the standard deviation $\sigma_n$ are asymptotically linear in $n$. 
\end{theorem}

To try to  apply Theorem~\ref{ass} to \eqref{Fap0} we set 
\[A(z,x)=0,\quad
 B(z,x)=(1-\a z)^{-\g/\a}, \quad C(z,x)=1+\frac{\b'x}\b\left(1-(1-\a z)^{-\b/\a}\right)\] 
 and  $\kappa=1+\g'/\b'$.
 Analytic perturbation condition of that theorem holds with any $r<1/\a$ 
 and 
 \[\rho=\frac1\a\left(1-\left(1+\frac{\b}{\b'}\right)^{-\a/\b}\right).\]
 (We may set $r=(1-\eps)/\a$ where $0<\eps<(1+\b/\b')^{-\a/\b}$, recall that the parameters are assumed to be non--negative.)  
 
 For the non--degeneracy 
 we find
 that if $\b'>0$ then
 \[\partial_zC(\rho,1)= -\b'( 1+ \b/\b')^{1+\a/\b},\quad \partial_xC(\rho,1)=-1.\]
 so that the condition holds and 
 \be\label{rhox}\rho(x)=\frac1\a\left(1-\left(1+\frac{\b}{\b'x}\right)^{-\a/\b}\right).\ee
  However,  if $\b'=0$ then $C(z,x)\equiv1$ and the non--degeneracy condition fails in that case. Thus, we consider the two cases separately.
 
 \section{The case $\b'>0$} 
 When $\b'>0$ we have 
 \begin{theorem}\label{no0s}
Let $\b'>0$ and $\a,\b\ge0$. Then: 
\begin{itemize}
\item[(i)] If $\a,\b>0$ then the  random variables $X_n$ satisfy
\be\label{clt}\frac{X_n-\mu_n}
{\sigma_n}
\stackrel d\rightarrow \N(0,1),\quad\mbox{as}\quad n\to\infty,\ee
with 
$\mu_n\sim-\mathfrak m(\rho)n$, $\sigma_n^2\sim-\mathfrak v(\rho)n$ where 
\begin{align*}\mathfrak m(\rho)&=-\frac{\a(1+\b/\b')^{-1-\a/\b}}{\b'(1-(1+\b/\b')^{-\a/\b})},\\
\mathfrak v(\rho)&=
 \frac{(1+\b/\b')^{-2-\a/\b}}{(\rho(1)\b')^2}\left(\frac{\b'}\a\left(1-(1+\b/\b')^{-\a/\b}\right)-1\right).\end{align*}
\item[(ii)] If $\a=0$, $\b>0$ then \eqref{clt} holds with 
\[\mathfrak m(\rho)=\frac{-\b}{(\b'+\b)\log(1+\b/\b')},\quad
\mathfrak v(\rho)=
\frac{\left(\frac{\b'}{\b}\log(1+\b/\b')-1\right)}{(\rho(1)(\b'+\b))^2}.
\]
\item[(iii)] If $\a>0$, $\b=0$  then the \eqref{clt} holds with 
\[\mathfrak m(\rho)=\frac{\a e^{-\a/\b'}}{\b'(e^{-\a/\b'}-1)},\quad \mathfrak v(\rho)=\frac{e^{-\a/\b'}}{(\rho(1)\b')^2}
\left(\frac{\b'}\a\left(1-e^{-\a/\b'}\right)-1\right).
\]
\end{itemize}
\end{theorem}
\proof
In view of the earlier discussion, to apply Theorem~\ref{ass} it suffices to verify the variability condition.

Let $\a,\b>0$ for now. 
 Note that  for a function $f$ and a constant $c\ne0$ we have $\mathfrak v(1/f)=-\mathfrak v(f)$ and  $\mathfrak v(cf)=\mathfrak v(f)$. Hence  $\mathfrak v(\rho(1)/\rho(x))=-\mathfrak v(\rho(x))$. After some calculations we get for a function $\rho(x)$ given by \eqref{rhox} that 
 \be\label{v}
 \mathfrak v(\rho)=
 \frac{(1+\b/\b')^{-2-\a/\b}}{(\rho(1)\b')^2}\left(\frac{\b'}\a\left(1-(1+\b/\b')^{-\a/\b}\right)-1\right).
 \ee
  The fraction is  positive and the parenthesized  quantity  is negative since
 \[\frac{\b'}\a\left(1-(1+\b/\b')^{-\a/\b}\right)-1=\frac{\b'\b}{\b\a}\left(1-(1+\b/\b')^{-\a/\b}\right)-1
 \] 
 and for $u,v>0$ we have
 \[\frac{1-(1+u)^{-v}}{uv}<1.\]
 (The last inequality is equivalent to $(1+u)^{-v}+uv-1>0$, and if for $v>0$ we let $g_v(u)$ to be the left hand side, we see that $g_v(0)=0$ and that $g_v$ is striclty incerasing for $u>0$.)
 Thus, the variability condition holds.

 It follows from Theorem~\ref{ass} that the corresponding random variables $X_n$ are asymptotically normal. 
 The variance of $X_n$ is  asymptotic to $-\mathfrak v(\rho)n$ with $\mathfrak v(\rho)$ given by \eqref{v}.  
 The expected value of $X_n$ is  asymptotic to 
 $-\mathfrak m(\rho)n$
 where $\mathfrak m(f)=f'(1)/f(1)$ and we similarly have $\mathfrak m(f(1)/f(x))=-\mathfrak m(f(x))$. 
 In our case
 \be\label{m}\mathfrak m(\rho)=-\frac{\a(1+\b/\b')^{-1-\a/\b}}{\b'(1-(1+\b/\b')^{-\a/\b})}.\ee
 (Note that $\mathfrak m(\rho)<0$ as it should be
  since
$1-(1+\b/\b')^{-\a/\b}>0$.)
This, proves the first part. 

 The remaining two parts follow by repeating the same calculations as for the case $\a,\b>0$ or by taking the limits as  the relevant parameters go to zero. Specifically, for the second part 
 letting $\a\to0$ in \eqref{v} and using
\[\lim_{\a\to0}\frac{1-(1+\b/\b')^{-\a/\b}}\a=\frac1\b\log(1+\b/\b')\]
yields
\[\mathfrak v(\rho)=
\frac{1}{(\rho(1)(\b'+\b))^2}\left(\frac{\b'}{\b}\log(1+\b/\b')-1\right)<0
\]
 since $\log(1+u)/u<1$ for $u>0$. Thus, the variability holds. 
 Upon taking $\a\to0$ in \eqref{m} we obtain
\[ \mathfrak m(\rho)=\frac{-\b}{(\b'+\b)\log(1+\b/\b')}\]
as claimed.  
Alternatively, with $\a\to0$, \eqref{Fap0} becomes
\be\label{a0}F(z,x)=\frac{e^{\g z}}
{\left(1+\frac{\b'x}\b\left(1-e^{\b z}\right)\right)^{1+\g'/\b'}}
\ee
and one can apply Theorem~\ref{ass} with 
\[A(z,x)=0,\quad
 B(z,x)=e^{\g z}, \quad C(z,x)=1+\frac{\b'x}\b\left(1-e^{\b z}\right).\] 
 For the last part, 
 as $\b\to0$ \eqref{v} becomes
\[\mathfrak v(\rho)=\frac{e^{-\a/\b'}}{(\rho(1)\b')^2}
\left(\frac{\b'}\a\left(1-e^{-\a/\b'}\right)-1\right)<0
\]
since $(1-e^{-u})/u<1$ for $u>0$. Thus, the variability holds. 
Similarly, letting $\b\to0$ in \eqref{m} yields
\[\mathfrak m(\rho)=\frac{\a e^{-\a/\b'}}{\b'(e^{-\a/\b'}-1)}\]
as desired.
\fpf

 The last case when $\b'>0$ is $\a=\b=0$.   
    Letting $\b\to0$ in \eqref{a0} gives
\be\label{f_3.3} F(z,x)=\frac{e^{\g z}}
{\left(1-\b'xz\right)^{1+\g'/\b'}}.
\ee
Here, $C(x)=1-\b'xz$ and $\rho(x)=(\b'x)^{-1}$. It follows that $\mathfrak v(\rho)=0$ and we find  that the variability condition 
fails.  Thus, Theorem~\ref{ass} is not applicable. 
However, we can identify the distribution of $X_n$ exactly. Before stating the result recall that if $r>0$ and $0<p<1$ than a random variable $Y$ has a negative binomial distribution with parameters $r$ and $p$ (abbreviated to $\NB(r,p)$) if 
\[\P(Y=k)= \frac{\Gamma(k+r)}{k!\Gamma(r)}(1-p)^kp^r={r+k-1\choose k}(1-p)^kp^r,\quad k=0,1,2,\dots
\]
We should note that some sources reserve the name negative binomial to  $r\in\mathbb N$ in which case $Y$ has the  interpretation as the number of failures before the $r$th success in the series of independent Bernoulli trials with parameter $p$.   Wikipedia refers to the $r>0$ case  as  P\'olya distribution (without giving a source). However, \cite{jk}  uses the  name negative binomial for any $r>0$ referring to the $r\in\mathbb N$ case  as Pascal distribution.  

We can now make the following statement.
\begin{theorem}
Let $\a=\b=0$. If $\g>0$ then the random variable $X_n$ has a $\NB(1+\g'/\b',\g/(\b'+\g))$ 
distribution, 
conditioned on its sum with an independent Poisson random variable $\Pois(\g/(\b'+\g))$ being equal to $n$. 
For the mean and the variance of $X_n$ we have 
\[\mu_n=n-\frac\g{\b'}(1+o(1)),\quad\sigma_n^2= \frac\g{\b'}\left(1+o(1)\right).
\]
If $\g=0$ then $X_n=n$ almost surely.\end{theorem}
\proof Consider $F(z,x)$ given by \eqref{f_3.3}. 
 Using the binomial expansion
\[(1-v)^{-\xi}=\sum_{k\ge0}{k+\xi-1\choose k}v^k
\] 
we find that 
\begin{align*}F(z,x)&=
\sum_{j\ge0}\frac{(\g z)^j}{j!}\sum_{k\ge0}{k+\g'/\b'\choose k}(\b'xz)^k\\&=\sum_{n\ge0}z^n\sum_{k=0}^n{k+\g'/\b'\choose k}(\b'x)^k\frac{\g^{n-k}}{(n-k)!}.
\end{align*}
Therefore, the PDF of the corresponding  random variable $X_n$ is
\[\E\ x^{X_n}=\frac{[z^n]F(z,x)}{[z^n]F(z,1)}=\sum_{k=0}^nx^kp_{n,k},\]
where
\[p_{n,k}=\frac{\b'^k{k+\g'/\b'\choose k}\frac{\g^{n-k}}{(n-k)!}}{\sum_{j=0}^n\b'^j{j+\g'/\b'\choose j}\frac{\g^{n-j}}{(n-j)!}},\quad0\le k\le n.
\]
When $\g=0$, $p_{n,k}=\delta_{k=n}$ so $X_n$  is a degenerate random variable equal to $n$.  Otherwise, 
re--writing $p_{n,k}$ as
\[p_{n,k}=
\frac{\left(\frac{\b'}{\b'+\g}\right)^k\left(\frac{\g}{\b'+\g}\right)^{1+\g'/\b'}
{k+\g'/\b'\choose k}
\frac{
e^{-\frac{\g}{\b'+\g}}}
{(n-k)!}\left(\frac{\g}{\b'+\g}\right)^{n-k}}
{\sum_{j=0}^n\left(\frac{\b'}{\b'+\g}\right)^j\left(\frac{\g}{\b'+\g}\right)^{1+\g'/\b'}
{j+\g'/\b'\choose j}
\frac{
e^{-\frac{\g}{\b'+\g}}}
{(n-j)!}\left(\frac{\g}{\b'+\g}\right)^{n-j}}
\]
and realizing that 
\[\left(\frac{\b'}{\b'+\g}\right)^k\left(\frac{\g}{\b'+\g}\right)^{1+\g'/\b'}
{k+\g'/\b'\choose k}
\]
 is  the  probability that a $\NB(1+\g'/\b',\g/(\b'+\g))$ random variable is equal to $k$, we conclude that  $X_n$ has a distribution of such a  random variable  conditioned on  its sum with an independent $\Pois(\g/(\b'+\g))$ being equal to $n$.

  To derive the numerical characteristics of this $X_n$ observe that
  \begin{align*}\E(n-X_n)
 & =\frac{\sum_{k=0}^n(n-k)\left(\frac{\b'}{\b'+\g}\right)^k\left(\frac{\g}{\b'+\g}\right)^{1+\g'/\b'}
{k+\g'/\b'\choose k}
\frac{
e^{-\frac{\g}{\b'+\g}}}
{(n-k)!}\left(\frac{\g}{\b'+\g}\right)^{n-k}}
{\sum_{j=0}^n\left(\frac{\b'}{\b'+\g}\right)^j\left(\frac{\g}{\b'+\g}\right)^{1+\g'/\b'}
{j+\g'/\b'\choose j}
\frac{
e^{-\frac{\g}{\b'+\g}}}
{(n-j)!}\left(\frac{\g}{\b'+\g}\right)^{n-j}}
\\& =\left(\frac{\g}{\b'+\g}\right)\frac{\sum_{k=0}^{n-1}\left(\frac{\b'}{\b'+\g}\right)^k
{k+\g'/\b'\choose k}
\frac1
{(n-1-k)!}\left(\frac{\g}{\b'+\g}\right)^{n-1-k}}
{\sum_{j=0}^n\left(\frac{\b'}{\b'+\g}\right)^j
{j+\g'/\b'\choose j}
\frac1
{(n-j)!}\left(\frac{\g}{\b'+\g}\right)^{n-j}}
\\&=\left(\frac{\g}{\b'+\g}\right)\frac{S_{n-1}}{S_n},
  \end{align*}
  where we have set
 \be\label{sn} S_n
=
\sum_{k=0}^{n}\left(\frac{\b'}{\b'+\g}\right)^k
\frac{\Gamma(k+1+\g'/\b')}
{k!(n-k)!}\left(\frac{\g}{\b'+\g}\right)^{n-k}
\ee
after using identity
\[{k+\g'/\b'\choose k}=\frac{\Gamma(k+1+\g'/\b')}{k!\Gamma(1+\g'/\b')}.\]
 Similarly,
\[\E(n-X_n)(n-1-X_n)=\left(\frac{\g}{\b'+\g}\right)^2\frac{S_{n-2}}{S_n}.
\]
It follows that 
\[\E X_n=n-\left(\frac{\g}{\b'+\g}\right)\frac{S_{n-1}}{S_n}\]
and that 
\begin{align*}\var(X_n)&=\var(n-X_n)=\E(n-X_n)(n-1-X_n)+\E(n-X_n)-\left(\E(n-X_n)\right)^2
\\&=\left(\frac{\g}{\b'+\g}\right)^2\frac{S_{n-2}}{S_n}+\left(\frac{\g}{\b'+\g}\right)\frac{S_{n-1}}{S_n}-\left(\frac{\g}{\b'+\g}\right)^2\frac{S^2_{n-1}}{S^2_n}
\\&=\left(\frac{\g}{\b'+\g}\right)\frac{S_{n-1}}{S_n}\left(1+\left(\frac{\g}{\b'+\g}\right)\left(\frac{S_{n-2}}{S_{n-1}}-\frac{S_{n-1}}{S_n}\right)\right).
\end{align*}
Using $\Gamma(k+1+\g'/\b')=(k+\g'/\b')\Gamma(k+\g'/b')$ we see from \eqref{sn} that 
\[S_n=\left(\frac{\g}{\b'+\g}\right)^n\frac1{n!}+\frac{\b'}{\b'+\g}S_{n-1}+\frac{\g'}{\b'}\sum_{k=1}^n\left(\frac{\b'}{\b'+\g}\right)^k\frac{\Gamma(k+\g'/\b')}{k!(n-k)!}\left(\frac{\g}{\b'+\g}\right)^{n-k}.
\]
We will show that the right--hand side is $\frac{\b'}{\b'+\g}S_{n-1}(1+o(1))$. It then follows that $S_{n-1}/S_n\to(\b'+\g)/\b' $ as $n\to\infty$ and thus
\[
\E X_n=n-\frac\g{\b'}(1+o(1)),\quad\var(X_n)= \frac\g{\b'}\left(1+o(1)\right)
\]
as desired.

For the first summand we trivially have
\[\left(\frac{\g}{\b'+\g}\right)^n\frac1{n!}=\frac{\g}{(\b'+\g)n}\left(\frac{\g}{\b'+\g}\right)^{n-1}\frac1{(n-1)!}
\le \frac{\g}{(\b'+\g)n}S_{n-1}.
\]
To show that 
\[\sum_{k=1}^n\left(\frac{\b'}{\b'+\g}\right)^k\frac{\Gamma(k+\g'/\b')}{k!(n-k)!}\left(\frac{\g}{\b'+\g}\right)^{n-k}=o(S_{n-1})\]
we split the sum as
\[\left(\sum_{k=1}^{k_n}\, +\sum_{k=k_n+1}^n\right)\left(\frac{\b'}{\b'+\g}\right)^k\frac{\Gamma(k+\g'/\b')}{k!(n-k)!}\left(\frac{\g}{\b'+\g}\right)^{n-k}
\]
with $k_n$ to be specified later. The second sum is bounded by
\begin{align*}&\frac{\b'}{(\b'+\g)(k_n+1)}\sum_{k=k_n+1}^n\left(\frac{\b'}{\b'+\g}\right)^{k-1}\frac{\Gamma(k+\g'/\b')}{(k-1)!(n-k)!}\left(\frac{\g}{\b'+\g}\right)^{n-k}
\\&\quad=\frac{\b'}{(\b'+\g)(k_n+1)}\sum_{k=k_n}^{n-1}\left(\frac{\b'}{\b'+\g}\right)^{k}\frac{\Gamma(k+1+\g'/\b')}{k!(n-1-k)!}\left(\frac{\g}{\b'+\g}\right)^{n-1-k}
\\&\quad \le \frac{\b'}{(\b'+\g)(k_n+1)}S_{n-1}
\end{align*}
which is $o(S_{n-1})$ as long as $k_n\to\infty$.

The first is bounded by
\begin{align*}&\frac{\b'}{\b'+\g}\sum_{k=1}^{k_n}\left(\frac{\b'}{\b'+\g}\right)^{k-1}\frac{\Gamma(k+\g'/\b')}{(k-1)!(n-k)!}\left(\frac{\g}{\b'+\g}\right)^{n-k}
\\&\quad=\frac{\b'}{\b'+\g}\sum_{k=0}^{k_n-1}\left(\frac{\b'}{\b'+\g}\right)^{k}\frac{\Gamma(k+1+\g'/\b')}{k!(n-1-k)!}\left(\frac{\g}{\b'+\g}\right)^{n-1-k}
\\&\quad=\frac{\b'}{\b'+\g}\left(\frac\g{\b'+\g}\right)^{n-k_n}\sum_{k=0}^{k_n-1}\left(\frac{\b'}{\b'+\g}\right)^{k}\frac{\Gamma(k+1+\g'/\b')}{k!(n-1-k)!}\left(\frac{\g}{\b'+\g}\right)^{k_n-1-k}
\end{align*}
Since $(n-1-k)!\ge (k_n-1-k)!(n-k_n)!$ the expression above is further bounded by
\[\frac{\b'}{\b'+\g}\left(\frac\g{\b'+\g}\right)^{n-k_n}\frac1{(n-k_n)!}S_{k_n-1}.
\]
Finally, since $S_n\ge \frac{\b'}{\b'+\g}S_{n-1}$ we obtain 
\[S_{k_n-1}=S_{n-1}\prod_{j=k_n}^{n-1} \frac{S_{j-1}}{S_j}\le\left(\frac{\b'+\g}{\b'}\right)^{n-k_n}S_{n-1}.
\]
Combining all of these  estimates shows that
\[\sum_{k=1}^{k_n}\left(\frac{\b'}{\b'+\g}\right)^k\frac{\Gamma(k+\g'/\b')}{k!(n-k)!}\left(\frac{\g}{\b'+\g}\right)^{n-k}\le
\frac{\b'(\g/\b')^{n-k_n}}{(\b'+\g)(n-k_n)!}
S_{n-1}=o(S_{n-1})
\] 
whenever $n-k_n\to\infty$. This completes the proof.
\fpf

\section{The case $\b'=0$}
We now consider the four cases corresponding to $\b'=0$. Letting $\b'\to0$ in \eqref{Fap0} yields
\be\label{Fbp0} F(z,x)=\frac{(1-\a z)^{-\g/\a}}
{\exp\left(\frac{\g'x}\b\left(1-(1-\a z)^{-\b/\a}\right)\right)}.
\ee
We begin with the simplest case, namely, $\b=0$.

\begin{theorem}\label{thm:b0} Assume that $\b=0$ in addition to $\b'=0$. 
If $\g'>0$ and 
\begin{itemize}
\item[(i)]\label{ane0}  $\a>0$ then 
\[\frac{X_n-(\g'/\a)\log n}{\sqrt{(\g'/\a)\log n}}\stackrel d\rightarrow \N(0,1).\] 
\item[(ii)]\label{thma0}  $\a=0$ then 
\[X_n\stackrel d= \Bin(n,\frac{\g'}{\g+\g'}).\]
\end{itemize}
If $\g'=0$ then $X_n$ is a degenerate random variable equal to $0$.
\end{theorem}
\proof 
Letting in \eqref{Fbp0} $\b\to0$ yields
\be\label{Fb0}F(z,x)=\frac{(1-\a z)^{-\g/\a}}
{\exp\left(\frac{\g'x}\a\log(1-\a z)\right)}=(1-\a z)^{-(\g'x+\g)/\a}.
\ee
If $\g'=0$ then $F(z,x)$ does not depend on $x$. Consequently,  $\E\ x^{X_n}=1$ and $X_n=0$ almost surely which proves the last assertion.

Otherwise, upon letting $\a\to0$, \eqref{Fb0} further simplifies to
\[
F(z,x)=e^{\g z+\g'xz}.
\]
The coefficient  $[z^n]F(z,x)$ is $(\g'x+\g)^n/n!$ which corresponds to $X_n$ having  the $\textrm{Bin}(n, \frac{\g'}{\g'+\g})$ distribution. This proves the second part. 


Finally, consider $\a>0$. 
Then, $F(z,x)$ given by \eqref{Fb0} falls within the variable exponent perturbation \cite[Theorem~IX.11,~p.~669]{fs} with 
\[A(z,x)=0,\quad B(z,x)=1,\quad C(z)=1-\a z,\]
and the exponent function  $\a(x)=(\g'x+\g)/\a$.
The variability condition $\a'(1)+\a''(1)\ne0$ of this theorem is satisfied  leading to asymptotic normality of $(X_n)$ with both $\mu_n$ and $\sigma_n^2$ asymptotic to $\frac{\g'}\a\log n$. 
\fpf

\begin{remark} We would like to comment that in the case when both $\b$ and $\b'$ are zero a  bit more precise information is readily available regardless of the value of $\a$ (or of $\a'$ as a matter of fact):  the recurrence \eqref{wpol} takes a simpler form
\[P_n(x)= ((\a'x+\a)n+\g'x+\g)P_{n-1}(x).\]
This means that
\[\E\ x^{X_n}=\frac{P_n(x)}{P_n(1)}=\prod_{k=1}^n\frac{(\a'x+\a)k+\g'x+\g}{(\a'+\a)k+\g'+\g}=\prod_{k=1}^n\frac{(\a'k+\g')x+\a k+\g}{(\a'+\a)k+\g'+\g}.
\] 
Thus,   $X_n$ is the sum of the independent indicators $I_k$, $1\le k\le n$, where 
\be\label{i}\P(I_k=1)=\frac{\a'k+\g'}{(\a'+\a)k+\g'+\g}.\ee
By elementary calculations
\[\E X_n=\frac{\a'}{\a'+\a}n+\frac{\g'\a-\g\a'}{(\a+\a')^2}\log n+O(1)\] 
and 
\[\var(X_n)=\frac{\a\a'}{(\a'+\a)^2}n+\frac{\a'^2\g+\a^2\g'-\a\a'(\g+\g')}{(\a+\a')^3}\log n+O(1).\]

In particular, when $\a'=0$ and $\a>0$ we recover the same asymptotic normality  as in Theorem~\ref{thm:b0} 
 while if both $\a$ and $\a'$ vanish then $X_n=\Bin(n,\frac{\g'}{\g+\g'})$ as is clear from \eqref{i}.  
We also note that  if
 $\a=0$ and $\a'>0$ then 
the mean is asymptotic to $n-\frac{\g}{\a'}\log n$ and the variance to $\frac{\g}{\a'}\log n$ so that in that case $(n-X_n-(\g/\a')\log n)/\sqrt{(\g/\a')\log n}$ is asymptotically $\N(0,1)$ as $n\to\infty$.
\end{remark}

It remains to consider the case $\b>0$. With $\a\to0$, \eqref{Fbp0} becomes
\[F(z,x)=\exp\left(\g z-\frac{\g'x}\b(1-e^{\b z})\right)=\exp\left\{\g z+\frac{\g'x}\b(e^{\b z}-1)\right\}.
\]
 This case was  treated in \cite{h} and gives the $\N(n/\log n,n/\log^2n)$ type of behavior. Specifically, we have
\begin
{theorem}[{\cite[Theorem~1]{h}}]
Suppose $\a=0$ and $\b>0$ in addition to $\b'=0$. Then, as $n\to\infty$, 
\[\frac{X_n-n/\log n}{\sqrt n/\log n}\stackrel d\rightarrow \N(0,1).\]
\end{theorem}

Finally, we consider the case when both $\a,\b>0$.
\begin{theorem}\label{abneo}
Assume $\a,\b>0$ in addition to $\b'=0$. Then
\[\frac{X_n-\mu_n}{\sigma_n}\stackrel d\rightarrow \N(0,1)\]
where 
\[\mu_n\sim\frac{(\g')^{\a/(\a+\b)}}\b(\a n)^{\b/(\a+\b)},\quad \sigma_n^2\sim
\frac{(\g')^{\a/(\a+\b)}\a (\a n)^{\b/(\a+\b)}}{\b(\a+\b)}.\]
\end{theorem}
\proof

Re--write \eqref{Fbp0} as 
\[F(z,x)=
\frac
{\exp\left(\frac{\g'x}\b\left((1-\a z)^{-\b/\a}\right)-1\right)}{(1-\a z)^{\g/\a}}.
\]
This has  the form discussed in \cite[Example~VIII.7]{fs} and  the asymptotics of the coefficients $[z^n]F(z,x)$  in a neighborhood of $x=1$  can be found by a version of the  saddle point method  developed by Wright in \cite{wr} and briefly discussed in \cite[Example~VIII.7 and Note~VIII.7]{fs}. Then a quasi--power theorem can be applied (see \cite[Theorem~IX.13]{fs} for a specific result applicable to this situation and \cite[Section~IX.8]{fs} for a broader discussion).  

Set
\be\label{fzx}f(z,x):=\frac{\g'x}\b\left((1-\a z)^{-\b/\a}-1\right)-\frac\g\a\log(1-\a z).\ee
The saddle point $r=r(x)$ is determined by the solution of 
\be\label{sad_pt_eq}(zf_z(z,x))_{z=r}=n
\ee
where the subscripts of $f$  denote its partial  derivatives  with respect to the indicated variables.
In the case at hand,  \eqref{sad_pt_eq}  becomes
\[\frac z{(1-\a z)^{1+\b/\a}}\left(\g'x+\g(1-\a z)^{\b/\a}\right)=n
\]
which leads, asymptotically, to
\be\label{ras}1-\a r(x)\sim\left(\frac{\g'x}{\a n}\right)^{\a/(\a+\b)}, \quad\mbox{i. e.}\quad r(x)\sim\frac1\a\left(1-\left(\frac{\g'x}{\a n}\right)^{\a/(\a+\b)}\right).\ee

From the general  theory it follows that the  PGF of a random variable $X_n$  encoded by $F(z,x)$ is given asymptotically by
 \[p_n(x)=\frac{F(r(x),x)}{F(r(1),1)}\left(\frac{r(1)}{r(x)}\right)^n(1+o(1)),
  \]
 where the $o(1)$ error is uniform over a small neighborhood  $\Omega$ of $x=1$.

 Let us set
   \[h_n(x):=f(r(x),x)-n\log r(x),\]
so that 
\[p_n(x)=\exp(h_n(x)-h_n(1))(1+o(1)).\]
 Then, by \cite[Theorem~IX.13]{fs} on generalized quasi--powers,  the corresponding random variables are asymptotically normal with mean asymptotic to $h'_n(1)$ and the variance asymptotic to $h_n'(1)+h_n''(1)$ provided 
 \be\label{h'''}\frac{h_n'''(x)}{(h'_n(1)+h''_n(1))^{3/2}}\rightarrow0,
 \ee
 uniformly over  a suitably small neighborhood $\Omega$ of $x=1$. Let us fixed such $\Omega$. To  verify \eqref{h'''} and evaluate the mean and the variance we need the first three derivatives of $h_n(x)$. First,  
   we get
 \beq*h_n'(x)&=&r'(x)f_z(z,x)_{z=r(x)}+f_x(z,x)_{z=r(x)}-n\frac{r'(x)}{r(x)}=f_x(r(x),x)
\eeq*
where the second equality follows from \eqref{sad_pt_eq} which implies that 
\[r'(x)\left(f_z(r(x),x)-\frac{n}{r(x)}\right)=0.
\]
Since $f(z,x)$ is a linear function of $x$,  $f_{xx}(z,x)=0$. Hence,   
\be\label{h''x}
h''_n(x)=r'(x)f_{zx}(z,x)_{z=r(x)}.
 \ee
Similarly,  $f_{xzx}(z,x)=0$ which  yields
\be\label{h'''x}h_n'''(x)=r''(x)f_{zx}(r(x),x)+(r'(x))^2f_{zzx}(r(x),x).\ee
For the function $f(z,x)$,  from \eqref{fzx} we see that 
\begin{align}\nonumber f_x(z,x)&=\frac{\g'}\b\left((1-\a z)^{-\b/\a}-1\right)\\
\nonumber f_{zx}(z,x)&=\g'(1-\a z)^{-1-\b/\a}
\\
\label{fzzx}f_{zzx}(z,x)&=\g'(\a+\b)(1-\a z)^{-2-\b/\a}.
\end{align}
Hence, using  \eqref{ras} (and writing $r=r(1)$, $h'_n=h'_n(1)$ etc.) we arrive at 
\be\label{h'}\E X_n\sim h'_n\sim\frac{\g'}\b\left((1-\a r)^{-\b/\a}-1\right)
\sim
\frac{(\g')^{\a/(\a+\b)}}\b(\a n)^{\b/(\a+\b)}.
\ee
For the variance we need to asymptotically evaluate $h_n''=r'f_{zx}(r,1)$ for which we will need the asymptotics of $r'$. To this end, observe that  
  implicit differentiation of \eqref{sad_pt_eq} 
  gives
 \[r'(x)f_z(r(x),x)+r(x)\Big(r'(x)f_{zz}(z,x)_{z=r(x)}+f_{zx}(z,x)_{z=r(x)}\Big)=0.
  \]
 which further gives, 
  \be\label{rhoprime}r'(x)=\frac{-r(x)f_{zx}(r(x),x)}{f_z(r(x),x)+r(x)f_{zz}(r(x),x)}.
  \ee
  Just as with earlier calculations,
 \be\label{fzx1}
  f_{zx}(r,1)\sim\g'(1-\a r)^{-1-\b/\a}\sim\g'\left(\left(\frac{\g'}{\a n}\right)^{\a/(\a+\b)}\right)^{-1-\b/\a}=\a n.
  \ee
  Further, from \eqref{fzx} we see that 
  \begin{align}\nonumber f_z(z,x)&=\g'x(1-\a z)^{-1-\b/\a}
  \\ \label{fzz} f_{zz}(z,x)&=\g'x(\a+\b)(1-\a z)^{-2-\b/\a}.
  \end{align}
  Consequently, using \eqref{ras} again we have
   \be\label{fzfzz}
 f_{z}(r,1)\sim\a n,\quad f_{zz}(r,1)\sim \frac{\b+\a}{(\g')^{\a/(\a+\b)}}(\a n)^{1+\a/(\a+\b)}.
\ee
Plugging these expressions in \eqref{rhoprime} and using $r\sim1/\a$ we obtain
\[r'\sim-\frac{n}{\a n+(\a+\b)\a^{-1}(\g')^{-\a/(\a+\b)}(\a n)^{1+\a/(\a+\b)}}\sim-\frac{(\g')^{\a/(\a+\b)}}{(\a+\b)(\a n)^{\a/(\a+\b)}}.
\]
Hence, from \eqref{h''x}
\[h_n''\sim-\frac{\a n(\g')^{\a/(\a+\b)}}{(\a+\b)(\a n)^{\a/(\a+\b)}}\sim-\frac{(\g')^{\a/(\a+\b)}}{\a+\b}(\a n)^{\b/(\a+\b)}.
\]
Combining this with \eqref{h'} we conclude that 
\[\var(X_n)
\sim(\g')^{\a/(\a+\b)}(\a n)^{\b/(\a+\b)}\left(\frac1\b-\frac1{\a+\b}\right)=\frac{(\g')^{\a/(\a+\b)}\a (\a n)^{\b/(\a+\b)}}{\b(\a+\b)}.
\]

It remains to verify \eqref{h'''} i.e. that $h_n'''(x)=o(n^{3\b/2(\a+\b)})$ uniformly over $\Omega$. From now on all approximations are understood to hold uniformly over $x\in\Omega$. We first note that while the earlier estimates were done for $x=1$, from the form of the expressions it is clear that they hold uniformly over $\Omega$ as well.  For example, it follows from \eqref{fzzx} and \eqref{ras} that 
\be\label{fzzxas}f_{zzx}(r(x),x)\sim\g'(\a+\b)(1-\a r(x))^{-2-\b/\a}\sim\g'(\a+\b)\left(\frac{\a n}{\g'x}\right)^{1+\a/(\a+\b)}
\ee
and that it holds uniformly over  $x\in\Omega$.
 
In the same fashion, by earlier estimates on $r'$ we find that  the second summand in the expression \eqref{h'''x} for $h_n'''(x)$  satisfies 
\be\label{2ndsumm}(r'(x))^2f_{zzx}(r(x),x)=O(n^{-2\a/(\a+\b)})\cdot O(n^{1+\a/(\a+\b)})=O(n^{\b/(\a+\b)})
\ee
uniformly.
It remains to bound $r''(x)f_{zx}(r(x),x)$. We know from \eqref{fzx1} that $f_{zx}(r(x),x)=O(n)$ so consider $r''(x)$.
 From \eqref{rhoprime} we see that 
\begin{align*}r''(x)&=-\frac{r'(x)f_{zx}(r(x),x)+r(x)r'(x)f_{zzx}(r(x),x)}{f_z(r(x),x)+r(x)f_{zz}(r(x),x)}
\\\quad&\qquad+\frac{r(x)f_{zx}(r(x),x)\frac d{dx}\left(f_{z}(r(x),x)+r(x)f_{zz}(r(x),x)\right)}{\left(f_z(r(x),x)+r(x)f_{zz}(r(x),x)\right)^2}.
\end{align*}
It follows from  \eqref{fzzxas}, \eqref{fzfzz} and the bound on $r'$ that the first term in absolute value is
\[O\left(\frac{n^{-\a/(\a+\b)}\left(n+n^{1+\a/(\a+\b)}\right)}{n+n^{1+\a/(\a+\b)}}\right)=O(n^{-\a/(\a+\b)}). \]
As for the second summand in the expression for $r''(x)$, first note that
\be\label{1stsum}\frac{r(x)f_{zx}(r(x),x)}{\left(f_z(r(x),x)+r(x)f_{zz}(r(x),x)\right)^2}=O\left(\frac n{n^{2+2\a/(\a+\b)}}\right)=O(n^{-1-2\a/(\a+\b)}).
\ee
To estimate the order of the derivative, we write
\begin{align*}\frac d{dx}\Big(f_{z}(r(x),x)&+r(x)f_{zz}(r(x),x)\Big)=2r'(x)f_{zz}(r(x),x)+f_{zx}(r(x),x)\\ &\quad+r(x)\Big(r'(x)f_{zzz}(r(x),x)+f_{zzx}(r(x),x)\Big).
\end{align*}
The new ingredient is $f_{zzz}(r(x),x)$ for which by \eqref{fzz} and \eqref{ras} we have
\[f_{zzz}(r(x),x)\sim\g'x(\a+\b)(2\a+\b)\left(\left(\frac{\g'x}{\a n}\right)^{\a/(\a+\b)}\right)^{-3-\b/\a}=O(n^{1+2\a/(\a+\b)}).
\]
Therefore, the derivative is of order
\[n^{-\a/(\a+\b)}n^{1+\a/(\a+\b)}+n+n^{-\a/(\a+\b)}n^{1+2\a/(\a+\b)}+n^{1+\a/(\a+\b)}=O(n^{1+\a/(\a+\b)}).
\]
Combining with \eqref{1stsum} we see that the second summand in $r''(x)$ is of order $n^{-\a/(\a+\b)}$. Since the first was of the same order we conclude that 
\[r''(x)=O(n^{-\a/(\a+\b)}).\]
Since by \eqref{fzfzz} $f_z(r(x),x)=O(n)$, in view of \eqref{2ndsumm} we see that each of the two terms comprising $h_n'''(x)$ in \eqref{h'''x} is of order $n^{\b/(\a+\b)}$. This completes verification of condition \eqref{h'''} and concludes the proof. 
\fpf

\section{Further remarks}

When $\a'\ne0$, the first characteristic equation in \eqref{char} is no longer separable, but  becomes the first order linear differential equation. In principle the underlying PDF can be solved, but usually does not give sufficiently explicit  form of the solution (see a comment in the Section~\ref{sec:mot}). Consequently,  one would have  to resort to treating solutions given in implicit form. However, when $\a=\b=0$ one still gets an explicit expression for the BGF which can then be analyzed.
 Namely,
\[F(z,x)=\frac{\exp\left\{\frac{\g}{\b'x}\left(1-\left(1-(\a'+\b')xz\right)^{\b'/(\a'+\b')}\right)\right\}}
{\left(1-(\a'+\b')xz\right)^{(\b'+\g')/(\a'+\b')}}.
\]
When $\g=0$ this is a degenerate random variable as is seen from the expansion of $(1-u)^{-\xi}$ in the powers of $u$. Otherwise, the PGF of $X_n$ is
\be\label{pgf_4.1}\frac{[z^n]F(z,x)}{[z^n]F(z,1)}=\exp\left\{\frac\g{\b'}\left(\frac1x-1\right)\right\}\frac{[z^n]H(z,x)}{[z^n]H(z,1)}\ee
where we have set
\[H(z,x):=
\frac{\exp\left\{-\frac{\g}{\b'x}\left(1-(\a'+\b')xz\right)^{\b'/(\a'+\b')}\right\}}{\left(1-(\a'+\b')xz\right)^{(\b'+\g')/(\a'+\b')}}.
\]
The first component on the right--hand side of \eqref{pgf_4.1} is the probability generating function of a negative of a $\Pois(\g/\b')$. When we expand $H$ as
\[
H(z,x)=\sum_{k\ge0}\left(-\frac\g{\b'x}\right)^k\frac{(1-(\a'+\b')xz)^{-\frac{(1-k)\b'+\g'}{\a'+\b'}}}{k!}.
\]
and then use the expansion of $(1-z)^{-\xi}$ 
we get
\begin{align*}[z^n]H(z,x)&=(\a'+\b')^nx^n\sum_{k\ge0}\left(-\frac\g{\b'x}\right)^k\frac{n^{\frac{(1-k)\b'+\g'}{\a'+\b'}-1}}{k!\Gamma(\frac{(1-k)\b'+\g'}{\a'+\b'})}\left(1+O(n^{-1})\right)\\&
=(\a'+\b')^n\sum_{k\ge0}x^{n-k}\frac{(-\g/{\b'})^k}{k!}\frac{n^{\frac{(1-k)\b'+\g'}{\a'+\b'}-1}}{\Gamma(\frac{(1-k)\b'+\g'}{\a'+\b'})}\left(1+O(n^{-1})\right)
\end{align*}
where the $O(n^{-1})$ can be replaced by a full expansion in the negative powers of $n$ as given  in \cite[Theorem~VI.1]{fs}  if needed. It follows that 
\[\left|\frac{[z^n]H(z,x)}{[z^n]H(z,1)}-x^n\right|=O(n^{-\b'/(\a'+\b')})\] 
uniformly over a small neighborhood of $x=1$. Hence $n-X_n$ is asymptotic to a $\Pois(\g/\b')$ random variable.

\end{document}